\documentclass[11pt]{article}
\usepackage{amsmath,amsxtra}
\usepackage{amssymb}
\usepackage{mathrsfs}
\usepackage{amsfonts}
\usepackage{amscd}
\usepackage[dvipdf]{graphics}
\usepackage[dvips]{graphicx}
\usepackage{graphicx,subfigure}
\usepackage{latexsym}
\usepackage{epstopdf}
\usepackage{latexsym}
\usepackage{amsmath}
\newcommand{\beqa}{\begin{eqnarray}}
\newcommand{\eeqa}{\end{eqnarray}}
\newcommand{\ba}{\begin{eqnarray*}}
\newcommand{\ea}{\end{eqnarray*}}

\date{}

\newtheorem{definition}{Definition}

\def\div{{\,\rm div \,}}

\def\R{\mathbb{R}}
\def\b{{\beta}}

\def\div{{\,\rm div \,}}

\def\be{{\begin{equation}}}
\def\ee{{\end{equation}}}
\def\Om{{\Omega}}
\def\na{{\nabla}}
\def\Ga{{\Gamma}}
\def\pl{{\partial}}
\def\qfq{{\quad\mbox{for}\quad}}
\def\beq{\arraycolsep=1.5pt\begin{eqnarray}}
\def\eeq{\end{eqnarray}}
\def\qiq{{\quad\mbox{in}\quad}}

\newfont{\Blackboard}{msbm10 scaled 1200}

\newfont{\roma}{cmr10 scaled 1200}

\def\qfq{{\quad\mbox{for}\quad}}
\def\<{{\langle}}
\def\>{{\rangle}}

\newtheorem{thm}{{}\hskip\parindent Theorem}[section]
\newtheorem{lem}{{}\hskip\parindent Lemma}[section]

\newtheorem{dfn}{{}\hskip\parindent Definition}[section]
\newtheorem{rem}{{}\hskip\parindent Remark}[section]

\def\be{\begin{equation}}
\def\ee{\end{equation}}
\def\beq{\arraycolsep=1.5pt\begin{eqnarray}}
\def\eeq{\end{eqnarray}}

\large
\title{\bf Energy decay and global solutions for a damped free boundary fluid-elastic structure interface model with variable coefficients in elasticity }
\date{}

\author{Yizhao Qin\quad
Peng-Fei Yao\thanks{Corresponding author.\ Email: pfyao@iss.ac.cn}\quad \\[0.3cm]
Key Laboratory of  Systems and Control\\
Institute of Systems Science,
Academy of Mathematics and Systems Science\\
Chinese Academy of Sciences, Beijing 100190, P. R. China\\
School of Mathematical Sciences\\
University of Chinese Academy of Sciences, Beijing 100049, China
}

\begin{document}
\maketitle
\footnote{This work is  supported by the National
Science Foundation of China, grants  no. 61473126 and no. 61573342, and Key Research Program of Frontier Sciences, CAS, no. QYZDJ-SSW-SYS011.}
\begin{quote}
\begin{small}
{\bf Abstract} \,\,\,We study a free boundary fluid-structure interaction model. In the model, a viscous incompressible fluid interacts with an elastic body via the common boundary. The motion of the fluid is governed by Navier-Stokes equations while the displacement of the elastic structure is described by variable coefficient wave equations. The dissipation is  placed on the common boundary between the fluid and the elastic body. Given small initial data, the global existence of the solutions of this system is proved and the exponential decay of solutions is obtained.
\\[3mm]
\\[3mm]
\end{small}
\end{quote}

\section{Introduction and Main Results}
\setcounter{equation}{0}
\hskip\parindent We consider a free boundary fluid-structure system which models the motion
of an elastic body moving and interacting with an incompressible viscous fluid
(see \cite{CS1,CS2,IKLT3, KT1,KT2}). This parabolic-hyperbolic system couples the Navier-Stokes equation
\be u_t-\Delta u + (u\cdot\na)u+\na p =0\quad\mbox{in}\quad \Om_f(t) \label{1.1}\ee
\be \na\cdot u =0\quad \mbox{in}\quad \Om_f(t)\label{1.2} \ee
with a wave equation with variable coefficients
\be w_{tt}-\div G(x)\na w +\b w =0\quad\mbox{in}\quad\Om_e\times(0,T),\label{1.3}\ee
where $\b>0$ is a constant. Without loss of generality, we set $\b=1$ from here on.
The Navier-Stokes equation is posed in the Eulerian framework and in a dynamic
domain $\Om_f(t),$  with $\Om_f(0)=\Om_f,$ while the wave equation is posed in the domain $\Om_e.$ The
geometry is such that $\pl\Om_e=\Ga_c$ is the common boundary of the domains, which is also the inner boundary of the domain $\Om_f,$ and we denote the outer boundary of $\Om_f$ as $\Ga_f=\pl\Om.$ This implies that $\pl\Om_f=\Ga_c\cup\Ga_f.$
Both domains $\Om_f$ and $\Om_e$ are assumed to be bounded and smooth (see \cite{CS1, IKLT3, KT1, KT2} for more
details).
In (\ref{1.3}),  $G(x)=(g_{ij}(x))_{3\times3}$ are symmetric, positive definite matrices for all $x\in\R^3,$ which are related to the elastic material.  For convenience, the entry $g_{ij}(x)$ is assumed to be smooth for $1\leq i,\,j\leq 3.$ The interaction is captured by
natural velocity and stress matching conditions on the free moving interface between the fluid and
the elastic body.

In the case of $G(x)=I,$ the identity matrix in $\R^{3\times 3},$ the existence of short time solutions was first established in \cite{CS1} and improved in \cite{IKLT1}  and  \cite{KT1} where there are no
dissipative  mechanisms on the interface. Then the global-in-time existence for the fluid-structure system with damping was established in \cite{IKLT2} and \cite{IKLT3} under the assumption of $G(x)=I$ for small initial data. In \cite{IKLT2}, the authors obtained the global-in-time existence in the following two cases: there are both internal damping and boundary dissipation or there is only boundary dissipation but with star-shaped condition for $\Om_e.$ Furthermore, in \cite{IKLT3}, the authors arrived at similar results by only imposing internal damping without any geometric conditions for domains.
For other topics on this system, see the short review in the introduction in \cite{IKLT3}.

Here we consider the fluid-structure system (\ref{1.1})-(\ref{1.3}) with the variabe coefficient elasticity structure described as the matrix $G(x)=\Big(g_{ij}(x)\Big)_{3\times3}$ to obtain
the global-in-time existence for small data.
The scheme, carried out in \cite{IKLT3}, and the geometrical in \cite{Yao 99} both play a key role here.
The geometrical approach was introduced in \cite{Yao 99} for the controllability of the wave equation with variable coefficients and extended in \cite{feng2001, LTY, TY, Yao 2007}, and many others, see \cite{Yao 2011}. The main advantages of this tool  are to provide great simplification (the Bochner technique, see \cite{wh 1989, GHL}) for the energy multipliers and  to yield the checkable geometrical assumptions for the variable coefficient problems by the curvature theory(refer to \cite{Yao 2011, wh 1989, dCa, GHL}), also see Remark \ref{r2.1} later.

To cope with the free boundary problem, we follow the approach in \cite{IKLT2}.
Let $\eta(x,t): \Omega\times(0,T)\rightarrow\Omega$ be the position function of the points in $\Om,$ which describes the different states of the system with respect to different time. With the help of position function, the incompressible Navier-Stokes equations can be rewritten as:
\begin{eqnarray}\label{eq2.1}
&&\left\{\begin{array}{lll}
\partial_{t}v^{i}-\partial_{j}(a^{jl}a^{kl}\partial_{k}v^i)+\partial_{k}(a^{ki}q)=0 \quad &\mbox{in}\quad \Omega_f\times (0,T),\\
a^{ki}\partial_{k}v^i=0  \quad &\mbox{in}\quad \Omega_f\times (0,T),
\end{array}\right. i=1,\,2,\,3,
\end{eqnarray}
where $v(x,t)$ and $q(x,t)$ denote the Lagrangian velocity and the pressure of the fluid over the initial domain $\Omega_f$, respectively.
This means that $v(x,t)=\eta_t(x,t)=u(\eta(x,t),t)$ and $q(x,t)=p(\eta(x,t),t)$ in $\Omega_f$. The matrix $\mathbf{a}(x,t)$ is defined as the inverse of the matrix $\nabla_x\eta(x,t)$, that is $\mathbf{a}(x,t)=(\nabla_x\eta(x,t))^{-1}$.
The wave equation for the displacement function $w(x,t)=\eta(x,t)-x$ is formulated in Lagrangian framework as
\begin{equation}\label{eq2.2}
w_{tt}^{i}-\div G\na w^i+w^i=0\quad\mbox{in}\quad \Omega_e\times (0,T),\quad i=1,2,3
\end{equation}
over the initial domain $\Om_e.$
We seek a solution $(v,w,q,\mathbf{a},\eta)$ to the  system \eqref{eq2.1}-\eqref{eq2.2}, where the matrix $\mathbf{a}=(a^{ij})(i,j=1,2,3)$ and $\eta\mid_{\Omega_f}$ are determined in the following way:
\begin{eqnarray}
&&\mathbf{a}_t=-\mathbf{a}:\nabla v:\mathbf{a},\quad\mbox{in}\quad \Omega_f\times(0,T), \label{eq2.3}\\
&&\eta_t=v \quad\mbox{in}\quad\Omega_f\times(0,T),\label{eq2.4}\\
&&\mathbf{a}(x,0)=I, \quad\eta(x,0)=x,\quad\mbox{in}\quad \Omega_f,  \label{eq2.5}
\end{eqnarray}
where the symbol $``:"$ stands for the usual multiplication between matrices.

Taking advantage of the matrix $\mathbf{a},$ we rewrite (\ref{eq2.1}) as
\begin{eqnarray*}
&&\left\{\begin{array}{lll}
\partial_{t}v-\div(\mathbf{a}:\mathbf{a}^{\mathbf{T}}:\nabla v)+\div(\mathbf{a}q)=0 \quad &\mbox{in}\quad \Omega_f\times (0,T),\\
tr(\mathbf{a}:\nabla v)=0  \quad &\mbox{in}\quad \Omega_f\times (0,T).\\
\end{array}\right.
\end{eqnarray*}
On the interface $\Gamma_c$ between $\Om_f$ and $\Om_e,$ we assume transmission boundary condition
\begin{equation}\label{eq2.6}
w_t=v-\gamma w_{\nu_\Lambda}\quad\mbox{on}\quad \Gamma_c\times(0,T),
\end{equation}
where the constant $\gamma>0$ and
$$w_{\nu_\Lambda}^{\mathbf{T}}=\nu^{\mathbf{T}}[G(x): \na w],$$
and the matching of stress
\begin{equation}\label{eq2.7}
w_{\nu_\Lambda}=(\mathbf{a}:\mathbf{a}^{\mathbf{T}}:\nabla v)\nu-q(\mathbf{a}\nu)\quad\mbox{on}\quad\Gamma_c\times(0,T),
\end{equation} where $\nu=(\nu_1,\nu_2,\nu_3)$ is the unit outward normal with respect to $\Om_e.$
On the outer fluid boundary $\Gamma_f,$ we impose the non-slip condition
\begin{equation}\label{eq2.8}
v=0\quad\mbox{on}\quad\Gamma_f\times(0,T).
\end{equation}
We supplement the system (\ref{eq2.1}) and (\ref{eq2.2})  with the initial conditions $v(x, 0)=v_0(x)$ and $(w(x, 0),w_t(x, 0))=(w_0(x),w_1(x))$
in $\Om_f$ and $\Om_e,$ respectively.
We employ the two classical  spaces
$$\mathbf{H}=\{u\in L^2(\Om_f): \div u=0,\quad u\cdot\nu|_{\Gamma_f}=0\}$$
and
$$\mathbf{V}=\{u\in H^1(\Om_f): \div u=0,\quad u|_{\Gamma_f}=0\}.$$
Based on the initial data $v_0$, the initial pressure $q_0$ is determined by solving the  problem
\begin{eqnarray}\label{eq2.9}
&&\left\{\begin{array}{lll}
\triangle q_0=-\partial_{i}v^{k}_{0}\partial_{k}v^{i}_{0}\quad &\mbox{in}\quad \Om_f,\\
\nabla q_{0}\cdot\nu=\triangle v_{0}\cdot\nu  \quad &\mbox{on}\quad \Gamma_f, \\
-q_0=-\partial_{j}v^{i}_{0}\nu_{j}\nu_{i}+w^{i}_{\nu_{\Lambda}}\nu_i,\quad &\mbox{on}\quad \Gamma_c. \\
\end{array}\right.
\end{eqnarray}
Let initial data $v_0\in\mathbf{V}\cap H^4(\Om_f),$ $w_0\in H^3(\Om_e),$ and $w_1\in H^2(\Om_e)$ be given. It is said that the {\it compatibility conditions} hold if $v_0,$ $w_0,$ and $w_1$ satisfy
\begin{eqnarray}
&&w_1=v_0-\gamma(w_0)_{\nu_{\Lambda}},\label{eq2.11}\\
&&w_{tt}(0)=v_t(0)-\gamma(w_1)_{\nu_{\Lambda}},\label{eq2.12} \\
&&w_{ttt}(0)=v_{tt}(0)-\gamma(w_{tt}(0))_{\nu_{\Lambda}},\label{eq2.13}
\end{eqnarray}
on $\Gamma_{c}$;
\begin{eqnarray}
&&(w_0)_{\nu_{\Lambda}}=\frac{\partial v_0}{\partial\nu},\label{eq2.14}\\
&&(w_1)_{\nu_{\Lambda}}=\frac{\partial v_t}{\partial\nu}(0)+\pl_t(\mathbf{a}:\mathbf{a}^{\mathbf{T}})(0):\nabla v_0\cdot\nu+\pl_t\mathbf{a}(0)\nu\cdot q(0),\label{eq2.15} \\
&&(w_{tt}(0))_{\nu_{\Lambda}}=\frac{\partial v_{tt}}{\partial\nu}(0)+2\pl_t(\mathbf{a}:\mathbf{a}^{\mathbf{T}})(0):\nabla v_t(0)\cdot\nu
+\pl_{tt}(\mathbf{a}:\mathbf{a}^{\mathbf{T}})(0):\nabla v_0\cdot\nu \nonumber\\
&&+2\pl_t\mathbf{a}(0)\nu\cdot q_t(0)+\pl_{tt}\mathbf{a}(0)\nu\cdot q(0),\label{eq2.16}
\end{eqnarray}
on $\Gamma_{c}$ as well; and
\begin{eqnarray}
&&v_0=0,\label{eq2.17}\\
&&\triangle v_0=\nabla q_0,\label{eq2.18} \\
&&-\partial_{t}a^{jk}(0)\partial_{jk}v_{0}^{i}-\partial_{j}(\partial_{t}a^{kj}(0)\partial_{k}v^{i}(0))
-\triangle\partial_{t}v^{i}(0)+\partial_{t}a^{ki}(0)\partial_{k}q(0)
+\partial_{it}q(0)=0\quad\quad\quad
\label{eq2.19}
\end{eqnarray}
on $\Gamma_f$.
Note that $w_{tt}(0),$ $w_{ttt}(0),$ $v_t(0)$ and $v_{tt}(0)$ can be represented by the initial data of the system with help of the structure of both Navier-Stokes and wave equations(see \cite{IKLT2} for detail).

For the existence of short time solutions, by following the treats in \cite{IKLT2}, we have the following.

\begin{thm}\label{th1.1}$(\cite{IKLT2})$
Assume that initial data $v_0\in\mathbf{V}\cap H^4(\Om_f),$ $w_0\in H^3(\Om_e),$ and $w_1\in H^2(\Om_e)$ are given small in some sense such that
the compatibility conditions (\ref{eq2.11})-(\ref{eq2.19}) are satisfied. Then, for any initial data $(v_0,w_0,w_1)$ satisfying the above conditions,
there exist a time $T_0>0$ and a unique solution $(v,w,q,\mathbf{a},\eta)$ to the system (\ref{1.1})-(\ref{1.3}) with the boundary conditions (\ref{eq2.6})-(\ref{eq2.8}), such that
\begin{eqnarray*}
&&v\in L^{\infty}([0,T_0];H^3(\Om_f)),\quad v_t\in L^{\infty}([0,T_0];H^2(\Om_f));\\
&&v_{tt}\in L^{\infty}([0,T_0];L^2(\Om_f)),\quad \nabla v_{tt}\in L^{2}([0,T_0];L^{2}(\Om_f));\\
&&\partial_{t}^{j}w\in\mathcal{C}([0,T_0];H^{3-j}(\Om_e)),\quad j=0,1,2,3,
\end{eqnarray*}
with $q\in L^{\infty}([0,T_0];H^2(\Om_f))$, $q_t\in L^{\infty}([0,T_0];H^1(\Om_f))$, $\mathbf{a},\mathbf{a}_t\in L^{\infty}([0,T_0];H^2(\Om_f))$, $\mathbf{a}_{tt}\in L^{\infty}([0,T_0];H^1(\Om_f))$, $\mathbf{a}_{ttt}\in L^{2}([0,T_0];L^{2}(\Om_f))$ and $\eta|_{\Om_f}\in\mathcal{C}([0,T_0];H^{3}(\Om_f))$.
\end{thm}

We omit the detail of the proof of Theorem \ref{th1.1}, but we still give the sketch of the construction of the short-time solutions following the procedure in \cite{IKLT2}. The local solutions of the fluid-structure interaction model are constructed by an iteration method. As shown in \cite{KT1,KT2}, it's sufficient to construct short-time solutions for a linear problem, which means that the problem with coefficients $\mathbf{a}(x,t)$ given and smooth, satisfies the postulated compatibility conditions. 
The solutions for the linear system is established as follows. First, we address the problem where $\mathbf{a}$ is smooth and independent of $t$ by a Galerkin procedure.
Then we obtain the existence of solutions with coefficients $\mathbf{a}(x,t)$ depending on $x$ and $t$ by contraction mapping principle. For the detail of the construction,
we refer to \cite{IKLT2}.

Next, we consider the hypothesis on the variable coefficients. We introduce
$$g=G^{-1}(x)\quad\mbox{for}\quad x\in\Om_e$$ as a Riemannian metric on $\Om_e$ and consider the couple $(\Om_e,g)$ as a Riemannian manifold. We denote by $g=\<\cdot,\cdot\>_g$ and $D$ the inner product and the covariant
differential of the metric $g,$ respectively, where
$$\<X,Y\>_g=\<G^{-1}(x)X,Y\>\qfq X,\,Y\in\R^3_x,\quad x\in\overline{\Om}_e,$$ where $\<\cdot,\cdot\>$ is the Euclidean metric of $\R^3.$   Note that $D$ is different from the Euclidean differential $\na$ unless $G(x)=I.$
\begin{dfn}
A vector field $H$ is said to be  {\it an \bf{escape vector field}} on $\overline{\Om}_e$ with respect to the metric $g$ (\cite{Yao 2011}), if it satisfies
\begin{equation*}
DH(X,X)\geqslant\varrho_{0}\mid X\mid_g^2,\quad\mbox{for any}\quad X\in\mathbb{R}^3_x,\quad x\in \overline{\Om}_e.
\end{equation*}
\end{dfn}

We need the following assumption.

{\bf (H)}\,\,There is an escape vector field $H$ in $\overline{\Om}_e$ with respect to the metric $g$ such that
\begin{equation*}
\langle H(x),\nu(x)\rangle\geqslant \gamma_0\quad\mbox{for}\quad x\in\overline{\Om}_e,
\end{equation*}
where $\gamma_0$ is a positive constant.

\begin{rem}\label{r2.1}If $G(x)=I,$ then $H=x-x_0$ for any $x_0\in\R^3$ given is an escape vector field. If $G$ is not the identity matrix, $H=x-x_0$ is no longer an escape vector field in general. The global existence of an escape vector field depends on
the sectional curvature of the metric $g.$  There does not exist such an escape vector field on $\overline{\Om}_e$ if $\Om_e$ contains a closed geodesic of the metric $g.$ For details, see $\cite[Section\, 3.2]{Yao 2011}.$
\end{rem}

Then we state the main result in this article. Our main result is the following.

\begin{thm}\label{th2.1}
Let the assumption ${\bf (H)}$ hold. Suppose that initial data $v_0\in\mathbf{V}\cap H^4(\Om_f),$ $w_0\in H^3(\Om_e),$ and $w_1\in H^2(\Om_e)$ are given small in the same sense as that in Theorem \ref{th1.1}, such that
the compatibility conditions (\ref{eq2.11})-(\ref{eq2.19}) hold.
Then, there exists a unique global smooth solution $(v,w,q,\mathbf{a},\eta)$ which satisfies
\begin{eqnarray*}
&&v\in L^{\infty}([0,\infty);H^3(\Om_f)),\quad v_t\in L^{\infty}([0,\infty);H^2(\Om_f));\\
&&v_{tt}\in L^{\infty}([0,\infty);L^2(\Om_f)),\quad \nabla v_{tt}\in L^{2}([0,\infty);L^{2}(\Om_f));\\
&&\partial_{t}^{j}w\in\mathcal{C}([0,\infty);H^{3-j}(\Om_e)),\quad j=0,1,2,3,
\end{eqnarray*}
with $q\in L^{\infty}([0,\infty);H^2(\Om_f))$, $q_t\in L^{\infty}([0,\infty);H^1(\Om_f))$, $\mathbf{a},\mathbf{a}_t\in L^{\infty}([0,\infty);H^2(\Om_f))$, $\mathbf{a}_{tt}\in L^{\infty}([0,\infty);H^1(\Om_f))$, $\mathbf{a}_{ttt}\in L^{2}([0,\infty);L^{2}(\Om_f))$ and $\eta|_{\Om_f}\in\mathcal{C}([0,\infty);H^{3}(\Om_f))$.
\end{thm}

\begin{rem}
 Given initial data small the total energy of the system decays exponentially.
\end{rem}

The exact meaning of the smallness of the initial data will be given later on.
The proof of Theorem \ref{th2.1} will be given in Section 4. Our paper is organized as following. In Section 2, we gather some lemmas and estimates which will be useful later on. In Section 3, we obtain the energy estimates and the superlinear estimates by geometric multiplier methods. The a-priori estimates attained in this section play a central role in our paper. Finally, in Section 4, the global existence and exponent decay of the energy of this fluid-structure system will be established.

\section{Preliminaries}
\setcounter{equation}{0}
\hskip\parindent We list some lemmas and estimates in the literature which are needed in the proof of Theorem \ref{th2.1}. Constant $C$ may be different  from line to line.

\begin{lem} \label{lem3.1}$\cite[Lemma\,\, 3.1]{IKLT1}$
Assume that $\parallel\nabla v\parallel_{L^{\infty}([0,T];H^2(\Om_f))}\leqslant M$. Let $p\in[1,\infty]$ and $1\leq i,\,j,\,k,\,l\leq 3.$
With $T\in[0,\frac{1}{CM}]$ where $C>0$ is large enough, the following statements hold:

$(i)$ $\parallel\nabla\eta\parallel_{H^2(\Omega_f)}\leqslant C\quad\mbox{for}\quad t\in[0,T];$

$(ii)$ $\parallel\mathbf{a}\parallel_{H^2(\Omega_f)}\leqslant C\quad\mbox{for}\quad t\in[0,T];$

$(iii)$ $\parallel\mathbf{a}_t\parallel_{L^p(\Omega_f)}\leqslant C\parallel\nabla v\parallel_{L^{p}(\Omega_f)}\quad\mbox{for}\quad t\in[0,T];$

$(iv)$ $\parallel\partial_{i}\mathbf{a}_t\parallel_{L^p(\Omega_f)}\leqslant C\parallel\nabla v\parallel_{L^{p_1}(\Omega_f)}\parallel\partial_{i}\mathbf{a}\parallel_{L^{p_2}(\Omega_f)}
+C\parallel\nabla\partial_{i}v\parallel_{L^p(\Omega_f)}, \quad t\in[0,T],$ where $1\leqslant p,p_1,p_2\leqslant\infty$
are given such that $\frac{1}{p}=\frac{1}{p_1}+\frac{1}{p_2}$;

$(v)$ $\parallel\partial_{ij}\mathbf{a}_t\parallel_{L^2(\Omega_f)}\leqslant C\parallel\nabla v\parallel_{H^1(\Omega_f))}^{\frac{1}{2}}\parallel\nabla v\parallel_{H^2(\Omega_f))}^{\frac{1}{2}}+C\parallel\nabla v\parallel_{H^2(\Omega_f)},\quad t\in[0,T];$

$(vi)$ $\parallel\mathbf{a}_{tt}\parallel_{L^2(\Omega_f)}\leqslant C\parallel\nabla v\parallel_{L^2(\Omega_f)}\parallel\nabla v\parallel_{L^{\infty}(\Omega_f)}+C\parallel\nabla v_t\parallel_{L^2(\Omega_f)} $ and    \\
$\parallel\mathbf{a}_{tt}\parallel_{L^3(\Omega_f)}\leqslant C\parallel v\parallel_{H^2(\Omega_f)}^2+C\parallel\nabla v_t\parallel_{L^3(\Omega_f)},$
for $t\in[0,T]$;

$(vii)$ $\parallel\mathbf{a}_{ttt}\parallel_{L^2(\Omega_f)}\leqslant C\parallel\nabla v\parallel_{H^1(\Omega_f)}^3+C\parallel\nabla v_t\parallel_{L^2(\Omega_f)}\parallel\nabla v\parallel_{L^{\infty}(\Omega_f)}+C\parallel\nabla v_{tt}\parallel_{L^2(\Omega_f),} t\in[0,T];$

$(viii)$ for every $\epsilon\in(0,\frac{1}{2}]$ and all $t\leqslant T^{*}=\min\{\frac{\epsilon}{CM^2},T\},$ we have
\begin{equation}\label{eq3.1}
\parallel\delta_{jk}-a^{jl}a^{kl}\parallel^2_{H^2(\Omega_f)}\leqslant \epsilon,\quad
\end{equation}
and
\begin{equation}\label{eq3.2}
\parallel\delta_{jk}-a^{jk}\parallel^2_{H^2(\Omega_f)}\leqslant \epsilon,\quad
\end{equation}
In particular, the form $a^{jl}a^{kl}\xi_{ij}\xi_{ik}$ satisfies the ellipticity estimates
\begin{equation}\label{eq3.3}
a^{jl}a^{kl}\xi_{ij}\xi_{ik}\geqslant \frac{1}{C}|\xi|^2,\quad \xi\in\mathbb{R}^{3\times3},
\end{equation}
for all $t\in[0,T^{*}]$ and $x\in\Omega_f$, provided $\epsilon\leqslant\frac{1}{C}$ with $C$ sufficiently large.
\end{lem}

From \cite{IKLT1}, we recall  pointwise a-priori estimates for variable coefficient Stokes system.
\begin{lem} \label{lem3.2}$\cite[Lemma\, 3.2]{IKLT1}$
Assume that $v$ and $q$ are solutions of the system
\begin{eqnarray*}
&&\left\{\begin{array}{lll}
\partial_{t}v^{i}-\partial_{j}(a^{jl}a^{kl}\partial_{k}v^i)+\partial_{k}(a^{ki}q)=0 \quad &\mbox{in}\quad \Omega_f\times (0,T),\\
a^{ki}\partial_{k}v^i=0  \quad &\mbox{in}\quad \Omega_f\times (0,T), \quad i=1,2,3\\
v^i=0,\quad\mbox{on}\quad \Gamma_f,\\
a^{jl}a^{kl}\partial_{k}v^{i}\nu_j-a^{ki}q\nu_k=(\omega^i)_{\nu_{\Lambda}},\quad\mbox{on}\quad \Gamma_c,
\end{array}\right.
\end{eqnarray*}
where $a^{ji}\in L^{\infty}(\Omega_f)$ satisfies lemma $\ref{lem3.1}$ with $\epsilon=\frac{1}{C}$ sufficiently small. Then we have
\begin{equation}\label{eq3.4}
\parallel v\parallel_{H^{s+2}(\Omega_f)}+\parallel q\parallel_{H^{s+1}(\Omega_f)}\leqslant C
\parallel v_t\parallel_{H^{s}(\Omega_f)}+C\parallel \omega_{\nu_{\Lambda}}\parallel_{H^{s+\frac{1}{2}}(\Gamma_c)},\quad s=0,1
\end{equation} for $t\in(0,T).$
Moreover, we  obtain
\begin{eqnarray}\label{eq3.5}
&&\parallel v_t\parallel_{H^2(\Omega_f)}+\parallel q_t\parallel_{H^1(\Omega_f)}\leqslant C\parallel v_{tt}\parallel_{L^2(\Omega_f)}+
C\parallel (w_t)_{\nu_{\Lambda}}\parallel_{H^\frac{1}{2}(\Gamma_c)}   \\
&&+C\parallel v\parallel_{H^2(\Omega_f)}^{\frac{1}{2}}\parallel v\parallel_{H^3(\Omega_f)}^{\frac{1}{2}}(\parallel v\parallel_{H^2(\Omega_f)}+\parallel q\parallel_{H^1(\Omega_f)})\quad\mbox{for }\quad t\in(0,T), \nonumber
\end{eqnarray}
where $T\leqslant\frac{1}{CM}$ and $C$ is sufficiently large.
\end{lem}

From now on, for simplicity,  we omit specifying the domains $\Omega_f$ and $\Omega_e$ in the norms involving the velocity $v$ and the displacement $w$. But we still emphasize the boundary domains $\Gamma_c$ and $\Gamma_f$.

Now we turn to the wave equations. Let $w$ be a solution to the wave system \eqref{eq2.2} satisfying the boundary condition
\eqref{eq2.6} on $\Gamma_c$. We write
\eqref{eq2.6}  as
\begin{equation*}
w_{\nu_{\Lambda}}=\frac{1}{\gamma}(v-w_t).
\end{equation*}
As in \cite{IKLT2}, by using the elliptic estimates and the Stokes type estimates in Lemma \ref{lem3.2}, we have the following estimates
\begin{eqnarray}\label{eq3.6}
\|w\|_{H^3}&\leq& C(\| w_{tt}\|_{H^1}
+\|w\|_{H^1}+\|v\|_{H^2}+\|w_t\|_{H^2})\quad\mbox{for}\quad t\in(0,T),
\end{eqnarray}
\begin{equation}\label{eq3.7}
\|w_t\|_{H^2}\leq C(\|w_{ttt}\|_{L^2}
+\parallel w_t\parallel_{L^2}+\parallel v_t\parallel_{H^1}+\parallel w_{tt}\parallel_{H^1})\qfq t\in(0,T),
\end{equation}
\begin{eqnarray}\label{eq3.8}
\parallel v\parallel_{H^3}+\parallel q\parallel_{H^2}
&\leqslant& C( \parallel v_t\parallel_{H^1}+\parallel v\parallel_{H^2}+\parallel w_{ttt}\parallel_{L^2}
+\parallel w_t\parallel_{L^2} \nonumber\\
&+&\parallel w_{tt}\parallel_{H^1})\qfq t\in(0,T),
\end{eqnarray}
\begin{equation}\label{eq3.9}
\parallel v\parallel_{H^2}+\parallel q\parallel_{H^1}\leqslant C(\parallel v_t\parallel_{L^2}+\parallel v\parallel_{H^1}+\parallel w_t\parallel_{H^1}),
\end{equation}
\begin{eqnarray}\label{eq3.10}
&&\parallel v_t\parallel_{H^2}+\parallel q_t\parallel_{H^1}
\leqslant C[\parallel v_{tt}\parallel_{L^2}+\parallel w_{ttt}\parallel_{L^2}
+\parallel w_t\parallel_{L^2}+\parallel v_t\parallel_{H^1} \nonumber\\
&&+\parallel w_{tt}\parallel_{H^1}+\parallel v\parallel_{H^3}^{\frac{1}{2}}(\parallel v_t\parallel_{L^2}+\parallel v\parallel_{H^1}+\parallel w_t\parallel_{H^1})^{\frac{3}{2}}],
\end{eqnarray}
for $t\in(0,T)$, where $T\leqslant\frac{1}{CM}$.

\section{A-priori estimates}
\setcounter{equation}{0}
\hskip\parindent
We derive a priori estimates for the global existence of  solutions to the dissipative fluid-structure system when the initial data are sufficiently small.

Suppose that
\begin{equation*}
\parallel v_0\parallel^2_{H^3}, \parallel v_t(0)\parallel^2_{H^1}, \parallel v_{tt}(0)\parallel^2_{L^2},
\parallel w_0\parallel^2_{H^3},\parallel w_1\parallel^2_{H^2}\leqslant \epsilon,
\end{equation*}
where $\epsilon>0$ is given small. This is the exact meaning of the smallness of the initial data we mention in Theorem \ref{th2.1} and Theorem \ref{th1.1}.

We need several auxiliary estimates involving different levels of energy as in the \cite{IKLT2}.
\subsection{First level estimates}
\hskip\parindent Let
\begin{equation}\label{eq4.1}
E(t)=\frac{1}{2}(\parallel v\parallel^2_{L^2}+\parallel w_t\parallel^2_{L^2}+\parallel w\parallel^2_{L^2}
+\|G^{1/2}\na w\parallel_{L^2}^2).
\end{equation}
$E(t)$ is the first level energy of the fluid-structure system.

As in \cite{IKLT2}, we have
\begin{lem}\label{l4.1}
Under the assumptions of Theorem \ref{th2.1}, the following inequality holds for $t\in[0,T]$
\begin{equation}\label{eq4.2}
E(t)+\int_0^{t}D(\tau)d\tau\leqslant E(0),
\end{equation}
where
\begin{equation}\label{eq4.3}
D(t)=\frac{1}{C}\parallel\nabla v(t)\parallel^2_{L^2}
+\gamma\parallel w_{\nu_{\Lambda}}(t)\parallel^2_{L^2(\Gamma_c)}
\end{equation}
is a dissipative term.
\end{lem}

To start with the application of the multiplier method, we need the two multiplier identities in  \cite[Theorems 2.1 and 2.2]{Yao 2011}. They are
\beq\label{eq4.4}
&&\div\{2\hat{H}(\hat{u})G(x)\nabla\hat{u}-(\mid\nabla_{g}\hat{u}\mid_g^2-\hat{u}_t^2)\hat{H}\}+2f\hat{H}(\hat{u}) \nonumber\\
&&=2[\hat{u}_{t}\hat{H}(\hat{u})]_t+2D\hat{H}(\nabla_{g}\hat{u},\nabla_{g}\hat{u})+(\hat{u}^2_t-\mid\nabla_{g}\hat{u}\mid_g^2)\div\hat{H},
\eeq
and
\beq\label{eq4.5}
&&\div[2p\hat{u}G(x)\nabla\hat{u}-\hat{u}^{2}G(x)\nabla p]+2fp\hat{u}\nonumber\\
&&=2p(\hat{u}\hat{u}_{t})_t+2p(\mid\nabla_{g}\hat{u}\mid_{g}^2-\hat{u}^2_t)-\hat{u}^{2}\div[G(x)\nabla p],
\eeq
where  $\hat{u}$ is a scalar function satisfying the wave equation $\hat{u}_{tt}=\div[G(x)\nabla\hat{u}]+f$ and $\hat{H}$ is a vector field while  $p$ is also a scalar function, where
$$\na_g\hat u=G(x)\na \hat u,\quad |\na_g\hat u|^2_g=\<G(x)\na\hat u,\na\hat u\>~~\mbox{and}~~\hat{H}(\hat{u})=\<\hat{H},\na\hat u\>.$$

We use $w$ as a multiplier as in \cite{IKLT2} to obtain

\begin{lem}\label{l4.2}
Suppose that the assumptions of Theorem \ref{th2.1} are true, it holds for $t\in[0,T]$ that
\beq\label{eq4.6}
&&\int_0^t\int_{\Omega_e}\< G(x)\na w,\nabla w\>dxd\tau+\int_0^t\parallel w\parallel^2_{L^2}d\tau\nonumber\\
&&\leq E(t)+E(0)+\int_0^t\parallel w_t\parallel^2_{L^2}d\tau+\int_0^t\langle w_{\nu_{\Lambda}},w\rangle_{L^2(\Gamma_c)}d\tau,
\eeq where $\< G(x)\na w,\na w\>=tr([G(x)\na w]^{\mathbf{T}}:\na w).$
\end{lem}

With help of  geometric multipliers, we have the following estimate which  plays the central role in deriving the total energy  of the solution to decay.

\begin{lem}\label{l4.3} Let the assumption ${\bf (H)}$ hold.
For any $\gamma>0$, there is a constant $C=C_{\gamma}>0$ such that
\begin{equation}\label{eq4.7}
E(t)+\int_0^{t} E(\tau)d\tau\leqslant CE(0),\quad\mbox{for}\quad t\in[0,T].
\end{equation}
In particular, $E(t)\leqslant CE(0)$.
\end{lem}

{\bf Proof} \,\,\,For simplicity, we further assume that the escape vector field $H$ satisfies
\begin{equation}\label{eq4.8}
DH(X,X)\geqslant\mid X\mid_g^2,\quad\mbox{for any}\quad X\in\mathbb{R}^n_x,\quad x\in \overline{\Omega}_e.
\end{equation}

First, we calculate in the components of $w$. In view of \eqref{eq2.2}, we take $H(w^i)$ as a multiplier and make use of  identity \eqref{eq4.4} by setting $\hat{u}=w^i$, $\hat{H}=H$ and $f=-w^i$
to have
\beq\label{eq4.9}
&&\div\{2H(w^i)G(x)\nabla w^i-(\mid\nabla_{g}w^i\mid_{g}^2-(w^i_t)^2)H\}-2w^{i}H(w^i)\nonumber\\
&&=
2[w^i_{t}H(w^i)]_t+2DH(\nabla_{g}w^i,\nabla_{g}w^i)+[(w^i_t)^2-\mid\nabla_{g}w^i\mid_g^2]\div H.
\eeq

Next, we integrate the above identities by parts over $\Om_e\times(0,t)$ for $i=1,$ $2,$ $3,$ and sum them to obtain
\begin{eqnarray}\label{eq4.10}
&&2\int_0^{t}\int_{\Omega_e}\sum_{i}DH(\nabla_{g}w^i,\nabla_{g}w^i)dxd\tau    \nonumber \\
&&+\int_0^{t}\int_{\Omega_e}(\mid w_t\mid^2-\< G(x)\na w,\nabla w\>-\mid w\mid^2) \div Hdxd\tau  \nonumber\\
&&+\int_0^{t}\int_{\Gamma_c}(\mid w\mid^2+\<G\na w,\na w\>)\langle H,\nu\rangle d\sigma d\tau
=-2\sum_{i}\int_{\Omega_e}w_t^{i}\cdot H(w^{i})dx\mid_0^t  \nonumber\\
&&+\int_0^{t}\int_{\Gamma_c}\mid w_t\mid^2\langle H,\nu\rangle d\sigma d\tau+2\sum_{i}\int^t_0\int_{\Gamma_c}H(w^{i})\cdot(w^i)_{\nu_{\Lambda}}d\sigma d\tau.
\end{eqnarray}
From \eqref{eq4.10}, we have
\begin{eqnarray}\label{eq4.11}
&&2\int_0^{t}\int_{\Omega_e}\<G\na w,\na w\>dxd\tau
+\int_0^{t}\int_{\Omega_e}(\mid w_t\mid^2-\<G\na w,\na w\>)-\mid w\mid^2)\div Hdxd\tau \nonumber\\
&&+\int_0^{t}\int_{\Gamma_c}(\mid w\mid^2+\<G\na w,\na w\>)\langle H,\nu\rangle d\sigma d\tau
\leqslant C(E(0)+E(t))+\int_0^{t}\int_{\Gamma_c}\mid w_t\mid^2\langle H,\nu\rangle d\sigma d\tau\nonumber\\
&&
+\epsilon_1\sum_i\int_0^{t}\int_{\Gamma_c}\mid H(w^i)\mid^{2}d\sigma d\tau
+C_{\epsilon_1}\int_0^{t}\int_{\Gamma_c}\mid w_{\nu_{\Lambda}}\mid^{2}d\sigma d\tau,
\end{eqnarray}
where $\epsilon_1>0$ is a sufficiently small constant to be determined.

We substitute $p=\frac{1}{2}\div H$ and $f=-w^i$ into \eqref{eq4.5} to have
\begin{eqnarray}\label{eq4.12}
&&\div\,[\div\,H w^{i}G(x)\nabla w^i-\frac{1}{2}(w^i)^{2}G(x)\nabla(\div\,H)]-\div\,H(w^i)^{2}
=\div\,H(w^i w^i_{t})_t \nonumber\\
&&+\div\,H[\mid\nabla_{g}w^i\mid_{g}^2-(w^i_t)^2]-\frac{1}{2}(w^{i})^{2}\div\,[G(x)\nabla(\div\,H)].
\end{eqnarray}
By a similar computation as in  \eqref{eq4.10}, we have
\begin{eqnarray}\label{eq4.13}
&&\int_0^{t}\int_{\Omega_e}(\mid w_t\mid^2-\<G\na w,\na w\>-\mid w\mid^2)\div\,Hdxd\tau
=\int_{\Omega_e}\div\,H\langle w,w_t\rangle dx\mid_{0}^{t}  \nonumber\\
&&-\frac{1}{2}\int_0^{t}\int_{\Omega_e}\mid w\mid^{2}\div\,[G(x)\nabla(\div\,H)]dxd\tau
-\int_0^{t}\int_{\Gamma_c}\div\,H\langle w,w_{\nu_{\Lambda}}\rangle d\sigma d\tau  \nonumber\\
&&+\frac{1}{2}\int_0^{t}\int_{\Gamma_c}\mid w\mid^{2}\langle\nabla_{g}(\div\,H),\nu\rangle_{g}d\sigma d\tau.
\end{eqnarray}
Moreover, we take $p=\frac{3}{2}$ and $f=-w^i$ in  identity \eqref{eq4.5} to get
\begin{equation}\label{eq4.14}
3\div\,[w^{i}G(x)\nabla w^i]-3(w^i)^{2}
=3(w^i w^i_{t})_t
+3[\mid\nabla_{g}w^i\mid_{g}^{2}-(w^i_t)^2].
\end{equation}
Integrating the above identities by parts leads to the following
\beq\label{eq4.15}
&&3\int_0^{t}\int_{\Omega_e}[\mid w_t\mid^2-\<G\na w,\na w\>-\mid w\mid^2]dxd\tau\nonumber\\
&&=3\int_{\Omega_e}\langle w,w_t\rangle dx\mid^t_0-3\int_0^{t}\int_{\Gamma_c}\langle w,w_{\nu_{\Lambda}}\rangle d\sigma d\tau.
\eeq

Insert  \eqref{eq4.13} into \eqref{eq4.11}, then add \eqref{eq4.15} together, and  we have
\begin{eqnarray}\label{eq4.16}
&&2\int_0^{t}\int_{\Omega_e}\<G\na w,\na w\>dxd\tau
+3\int_0^{t}\int_{\Omega_e}[\mid w_t\mid^2-\<G\na w,\na w\>-\mid w\mid^2]dxd\tau \nonumber\\
&&+\int_0^{t}\int_{\Gamma_c}(\mid w\mid^2+\<G\na w,\na w\>)\langle H,\nu\rangle d\sigma d\tau
\leqslant C(E(0)+E(t)) \nonumber\\
&&+\int_0^{t}\int_{\Gamma_c}\mid w_t\mid^2\langle H,\nu\rangle d\sigma d\tau
+\epsilon_1C\int_0^{t}\int_{\Gamma_c}|\na w|^2d\sigma d\tau  \nonumber\\
&&+C_{\epsilon_1,\epsilon_2}\int_0^{t}\int_{\Gamma_c}\mid w_{\nu_{\Lambda}}\mid^{2}d\sigma d\tau
+C\int_0^{t}\int_{\Omega_e}\mid w\mid^2dxd\tau  \nonumber\\
&&+\epsilon_2\int_0^{t}\int_{\Gamma_c}\mid w\mid^2d\sigma d\tau
+C\int_0^{t}\int_{\Gamma_c}\mid w\mid^2d\sigma d\tau.
\end{eqnarray}
Using the assumption $\langle H,\nu\rangle\geqslant \gamma_0$ in (\ref{eq4.16}) and letting $\epsilon_1,$ $\epsilon_2$ be sufficiently small such that
$$ C\epsilon_1+\epsilon_2\leq\gamma_0/2,    $$ we obtain
\begin{eqnarray}\label{eq4.17}
&&\frac{3}{2}\int_0^{t}\int_{\Omega_e}\<G\na w,\na w\>dxd\tau
+3\int_0^{t}\int_{\Omega_e}[\mid w_t\mid^2-\<G\na w,\na w\>-\mid w\mid^2]dxd\tau  \nonumber\\
&&+\frac{\gamma_0}{2}\int_0^{t}\int_{\Gamma_c}(\mid w\mid^2+\<G\na w,\na w\>)d\sigma d\tau
\leqslant C(E(0)+E(t))  \nonumber\\
&&+\int_0^{t}\int_{\Gamma_c}\mid w_t\mid^2\langle H,\nu\rangle d\sigma d\tau
+C\int_0^{t}\int_{\Gamma_c}\mid w_{\nu_{\Lambda}}\mid^{2}d\sigma d\tau  \nonumber\\
&&+C\int_0^{t}\int_{\Omega_e}\mid w\mid^2dxd\tau
+C\int_0^{t}\int_{\Gamma_c}\mid w\mid^2d\sigma d\tau.
\end{eqnarray}
Now, we multiply \eqref{eq4.6} by $3-2\epsilon_3$, where $\epsilon_3\in(0,1)$ is given small, to get
\begin{eqnarray}\label{eq4.18}
&&(3-2\epsilon_3)\int_0^t\int_{\Omega_e}\<G\na w,\na w\>dxd\tau
+(3-2\epsilon_3)\int_0^t\parallel w\parallel^2_{L^2}d\tau \nonumber\\
&&\leqslant (3-2\epsilon_3)(E(t)+E(0))
+(3-2\epsilon_3)\int_0^t\parallel w_t\parallel^2_{L^2}d\tau  \nonumber\\
&&+\frac{\gamma_0}{8}\int_0^t\int_{\Gamma_c}\mid w\mid^{2}d\sigma d\tau
+C\int_0^t\int_{\Gamma_c}\mid w_{\nu_{\Lambda}}\mid^{2}d\sigma d\tau.
\end{eqnarray}
Add \eqref{eq4.18} to \eqref{eq4.17}, and we have
\begin{eqnarray}\label{eq4.19}
&&(\frac{3}{2}-2\epsilon_3)\int_0^{t}\int_{\Omega_e}\<G\na w,\na w\>dxd\tau
+2\epsilon_3\int_0^{t}\int_{\Omega_e}\mid w_t\mid^{2}dxd\tau \nonumber\\
&&\leqslant C(E(0)+E(t))
+C\int_0^{t}\int_{\Omega_e}\mid w\mid^2dxd\tau
+\int_0^{t}\int_{\Gamma_c}\mid w_t\mid^2\langle H,\nu\rangle d\sigma d\tau  \nonumber\\
&&+C\int_0^{t}\int_{\Gamma_c}\mid w_{\nu_{\Lambda}}\mid^{2}d\sigma d\tau
+C\int_0^{t}\int_{\Gamma_c}\mid w\mid^2d\sigma d\tau.
\end{eqnarray}
Thus we have
\begin{eqnarray}\label{eq4.20}
&&\int_0^{t}\int_{\Omega_e}\<G\na w,\na w\>dxd\tau
+\int_0^{t}\int_{\Omega_e}\mid w\mid^2dxd\tau
+\int_0^{t}\int_{\Omega_e}\mid w_t\mid^{2}dxd\tau \nonumber\\
&&\leqslant C(E(0)+E(t))
+C\int_0^{t}\int_{\Omega_e}\mid w\mid^2dxd\tau  \nonumber\\
&&+C\int_0^{t}\int_{\Gamma_c}(\mid w_t\mid^2+\mid w_{\nu_{\Lambda}}\mid^{2}) d\sigma d\tau
+C\int_0^{t}\int_{\Gamma_c}\mid w\mid^2d\sigma d\tau.
\end{eqnarray}

To complete the proof of \eqref{eq4.7},  the lower order terms in (\ref{eq4.20}) have to be absorbed. This will be done by a compactness-uniqueness argument in the following lemma.
\begin{lem}\label{l3.4}
Under the assumptions of Lemma \ref{l4.3}, for any solution $(v,w)$ of the fluid-structure interaction system we deal with, we have that for all $t\in [0,T],$
\begin{eqnarray*}
\int_0^{t}\int_{\Omega_e}\mid w\mid^2dxd\tau+\int_0^{t}\int_{\Gamma_c}\mid w\mid^{2}d\sigma d\tau
\leqslant C(E(0)+E(t))+\int_0^{t}\int_{\Gamma_c}(\mid w_t\mid^2+\mid w_{\nu_{\Lambda}}\mid^{2}) d\sigma d\tau).
\end{eqnarray*}
\end{lem}

{\bf Proof of Lemma \ref{l3.4}}
We prove the estimate by contradiction. Suppose that there exists a solution sequence
 $(v_k,w_k)$ of the fluid-structure interaction system we study, such that
\beq
&&\int_0^{t}\int_{\Omega_e}\mid w_k\mid^2dxd\tau+\int_0^{t}\int_{\Gamma_c}\mid w_k\mid^{2}d\sigma d\tau=1, \label{3.22}\\
&&E(w_k,0)+E(w_k,t)+\int_0^{t}\int_{\Gamma_c}(\mid(w_k)_t\mid^2+\mid(w_k)_{\nu_{\Lambda}}\mid^{2}) d\sigma d\tau
<\frac{1}{k}\label{3.21}\eeq for $k=1,$ $2,$ $3,$ $\cdots,$
where $E(w_k,t)$ are obtained by replacing $w$ in $E(t)$ with $w_k.$

Taking $(v,w)=(v_k,w_k)$ in (\ref{eq4.20}), we deduce that the sequence
$\{w_k\}$ is bounded in $H^1(\Omega_e\times(0,t))$. Thus there is $\tilde{w}\in H^1(\Omega_e\times(0,t))$ such that
$$ w_k\rightharpoonup \tilde{w}\qiq H^1(\Omega_e\times(0,t)) \quad\mbox{weakly.}$$
Then $\tilde{w}$ is a weak solution to \eqref{1.3} in $\Omega_e\times(0,t)$ and
$$ w_k\rightarrow \tilde{w}\qiq L^2(\Omega_e\times(0,t)) \quad\mbox{strongly,}$$
$$\int_0^{t}\int_{\Gamma_c}\mid w_k\mid^{2}d\sigma d\tau\rightarrow\int_0^{t}\int_{\Gamma_c}\mid\tilde{w}\mid^{2}d\sigma d\tau,$$ as $k\rightarrow\infty.$
We obtain
\be
\int_0^{t}\int_{\Omega_e}\mid\tilde{w}\mid^2dxd\tau+\int_0^{t}\int_{\Gamma_c}\mid\tilde{w}\mid^{2}d\sigma d\tau=1.\label{3.23}
\ee

Moreover, it follows from (\ref{3.21}) that
$$\int_0^{t}\int_{\Gamma_c}\mid(w_k)_{\nu_{\Lambda}}\mid^{2}d\sigma d\tau\rightarrow 0,\quad
E(w_k,0)\rightarrow 0,\quad
\int_0^{t}\int_{\Gamma_c}\mid(w_k)_t\mid^2d\sigma d\tau\rightarrow 0$$ as $k\rightarrow\infty$.
Thus, we have
$$\tilde{w}_{\nu_{\Lambda}}\mid_{\Gamma_c}
=0,\quad \tilde{w}\mid_{\Gamma_c}=0,\quad \tilde{w}(0)=0.$$ As a consequence of \cite[Theorem 2.15 and Lemma 2.5]{Yao 2011}, we obtain
 $$\tilde{w}\equiv 0\qiq \Omega_e\times(0,t),$$ which contradicts (\ref{3.23}). \hfill$\Box$

After we obtain the Lemma \ref{l3.4}, we continue the proof of Lemma \ref{l4.3}.

Using \eqref{eq2.6} and Lemma \ref{l3.4} in (\ref{eq4.20}), we obtain
\begin{eqnarray}\label{eq4.23}
\int_0^{t}E(\tau)d\tau
\leqslant C(E(0)+E(t))+C\int_0^{t}\int_{\Omega_f}\mid\nabla v\mid^{2}dxd\tau
+C_{\gamma}\int_0^{t}\int_{\Gamma_c}\mid w_{\nu_{\Lambda}}\mid^{2}d\sigma d\tau.
\end{eqnarray}

Multiply \eqref{eq4.23} with $\epsilon_4>0$ small enough,  add it to \eqref{eq4.2}, and  then \eqref{eq4.7} follows. \hfill$\Box$

\begin{rem}
If the solution of this system can exist for all time and the matrix $\mathbf{a}$ can stay as close as possible to identity so that it is uniformly elliptic, then the above lemma implies that the first level energy $E(t)$ decays exponentially.
\end{rem}

\subsection{Second level estimates}
\hskip\parindent
Now we introduce the second level energy of the fluid-structure system
\begin{equation*}
E_1(t)=\frac{1}{2}(\parallel v_t\parallel^2_{L^2}+\parallel w_t\parallel^2_{L^2}+\parallel w_{tt}\parallel^2_{L^2}
+\|G^{1/2}\na w_t\|^2_{L^2}).
\end{equation*}
and the corresponding dissipation
\begin{equation*}
D_1(t)=\frac{1}{C}\parallel\nabla v_t(t)\parallel^2_{L^2}
+\gamma\parallel(w_t)_{\nu_{\Lambda}}\parallel^2_{L^2(\Gamma_c)}
\end{equation*}
To obtain similar estimates to (\ref{eq4.2}), (\ref{eq4.6}) and (\ref{eq4.7}) which are obtained in the above subsection, we differentiate the whole system in time following the approach of \cite{IKLT2}.

\begin{eqnarray}
&&v^{i}_{tt}-\partial_t\partial_{j}(a^{jl}a^{kl}\partial_{k}v^i)+\partial_t\partial_{k}(a^{ki}q)=0 \quad\mbox{in}\quad \Omega_f\times (0,T),\label{eq4.30}\\
&&a^{ki}\partial_{k}v^i_t+\partial_ta^{ki}\partial_{k}v^i=0  \quad\mbox{in}\quad \Omega_f\times (0,T), \quad i,j=1,2,3 \label{eq4.31}\\  \nonumber
\end{eqnarray}
and
\begin{equation}\label{eq4.32}
w_{ttt}^{i}-\div G(x)\na w^i_t+w^i_t=0\quad\mbox{in}\quad \Omega_e\times (0,T),\quad i=1,2,3.
\end{equation}

Moreover, the boundary conditions are changed into the form:
\begin{eqnarray}
&&w_{tt}=v_t-\gamma(w_t)_{\nu_\Lambda},\quad\mbox{on}\quad \Gamma_c\times(0,T),\label{eq4.33}\\
&&(w_t)_{\nu_\Lambda}^i=\partial_t(a^{jl}a^{kl}\partial_{k}v^{i})\nu_j-\partial_t(a^{ki}q)\nu_k,
\quad\mbox{on}\quad\Gamma_c\times(0,T),\label{eq4.34}\\  \nonumber
\end{eqnarray}
where $i,j,k,l,\alpha,\beta=1,2,3 $;

and
\begin{equation}\label{eq4.35}
v_t=0,\quad\mbox{on}\quad \Gamma_f\times(0,T).
\end{equation}

\begin{lem}\label{l4.5}
Under the assumptions of Theorem \ref{th2.1}, the following energy inequality holds for $t\in[0,T]$
\begin{equation}\label{eq4.36}
E_1(t)+\int_0^{t}D_1(\tau)d\tau\leqslant E_1(0)+\int^t_0(R_1(\tau),v_t(\tau))d\tau,
\end{equation}
where
\begin{eqnarray}
\int^t_0(R_1(\tau),v_t(\tau))d\tau=&-&\int_0^{t}\int_{\Omega_f}\partial_t(a^{jl}a^{kl})\partial_{k}v^i\partial_{j}v^{i}_{t}dxd\tau\\ \nonumber
&+&\int_0^{t}\int_{\Omega_f}\partial_ta^{ki}q\partial_kv^{i}_{t}dxd\tau
-\int_0^{t}\int_{\Omega_f}\partial_ta^{ki}\partial_tq\partial_kv^{i}dxd\tau.
\end{eqnarray}\label{eq4.37}
\end{lem}
{\bf Proof.} Take $L^2$ inner product with $v^i_t$ and $w^i_{tt}$ to \eqref{eq4.30} and \eqref{eq4.32}, respectively. Using the boundary conditions \eqref{eq4.33}-\eqref{eq4.35}, we
get that
\begin{eqnarray}\label{eq4.38}
\frac{1}{2}\frac{d}{dt}\parallel v_t\parallel^2_{L^2}+\int_{\Omega_f}a^{jl}a^{kl}\partial_{k}v^i_t\partial_{j}v^{i}_{t}dx
+\int_{\Omega_f}\partial_t(a^{jl}a^{kl})\partial_{k}v^i\partial_{j}v^{i}_{t}dx  \\ \nonumber
-\int_{\Omega_f}\partial_t(a^{ki}q)\partial_kv^{i}_{t}dx+\int_{\Gamma_c}\langle(w_t)_{\nu_\mathcal{A}},v_t\rangle d\sigma=0
\end{eqnarray}
and
\begin{eqnarray}\label{eq4.39}
\frac{1}{2}\frac{d}{dt}(\parallel w_{tt}\parallel^2_{L^2}+\parallel w_{t}\parallel^2_{L^2}+\int_{\Omega_e}\<G\na w_t,\na w_t\>dx)
-\int_{\Gamma_c}\langle(w_t)_{\nu_\Lambda},v_t-\gamma(w_t)_{\nu_\Lambda}\rangle d\sigma=0.\quad\quad
\end{eqnarray}
Add \eqref{eq4.38} and \eqref{eq4.39} together and integrate in time from $0$ to $t$. Thanks to the ellipticity of $\mathbf{a}(x,t)$,
we get
\begin{equation}\label{eq4.40}
E_1(t)+\int_0^{t}D_1(\tau)d\tau\leqslant E_1(0)-\int^t_0\int_{\Omega_f}\partial_t(a^{jl}a^{kl})\partial_{k}v^i\partial_{j}v^{i}_{t}dxd\tau
+\int^t_0\int_{\Omega_f}\partial_t(a^{ki}q)\partial_kv^{i}_{t}dxd\tau
\end{equation}

Taking advantage of \eqref{eq4.31}, we can rewrite \eqref{eq4.40} as
\begin{eqnarray}\label{eq4.41}
\int^t_0\int_{\Omega_f}\partial_t(a^{ki}q)\partial_kv^{i}_{t}dxd\tau
=\int^t_0\int_{\Omega_f}\partial_ta^{ki}q\partial_kv^{i}_{t}dxd\tau
+\int^t_0\int_{\Omega_f}a^{ki}\partial_{t}q\partial_kv^{i}_{t}dxd\tau  \nonumber\\
=\int^t_0\int_{\Omega_f}\partial_ta^{ki}q\partial_kv^{i}_{t}dxd\tau
-\int^t_0\int_{\Omega_f}\partial_ta^{ki}\partial_tq\partial_kv^{i}dxd\tau
\end{eqnarray}

We submit \eqref{eq4.41} into \eqref{eq4.40} and  we obtain (\ref{eq4.36}).  \hfill$\Box$\\

Next, we continue carrying out the multiplier method as before. First, as what we do in Lemma 3.2, we use $w_t$ as a
multiplier and get an analogy.
\begin{lem}\label{l4.6}
Under the hypotheses of Theorem \ref{th2.1}, we have for all $t\in[0,T]$ that
\beq\label{eq4.42}
&&\int_0^t\int_{\Omega_e}\<G\na w_t,\na w_t\>dxd\tau+\int_0^t\parallel w_t\parallel^2_{L^2}d\tau\nonumber\\
&&\leqslant
E_1(t)+E_1(0)+\int_0^t\parallel w_{tt}\parallel^2_{L^2}d\tau
+\int_0^t\langle(w_t)_{\nu_{\Lambda}},w_t\rangle_{L^2(\Gamma_c)}d\tau.
\eeq
\end{lem}

We go on applying the multiplier method to derive an analogy to Lemma 3.3. This time, we regard $H(w^i_t)$ and $\div\,Hw^i_t$ as multipliers. Through a similar procedure, by Lemma 3.5, Lemma 3.6 and boundary condition \eqref{eq4.33}, we obtain the estimates in Lemma 3.7.
\begin{lem}\label{l4.7}
For all $t\in[0,T]$, it holds under the assumptions of Theorem \ref{th2.1} that
\begin{equation}\label{eq4.43}
E_1(t)+\int_0^tE_1(\tau)d\tau\leqslant
C(E_1(0)+\int^t_0(R_1(\tau),v_t(\tau))d\tau).
\end{equation}
\end{lem}

\subsection{Third level estimates}
\hskip\parindent
We introduce the third level energy
\begin{equation*}
E_2(t)=\frac{1}{2}(\parallel v_{tt}\parallel^2_{L^2}+\parallel w_{tt}\parallel^2_{L^2}+\parallel w_{ttt}\parallel^2_{L^2}
+\|G^{1/2}\na w_{tt}\|^2_{L^2}).
\end{equation*}
and the corresponding dissipation
\begin{equation*}
D_2(t)=\frac{1}{C}\parallel\nabla v_{tt}(t)\parallel^2_{L^2}
+\gamma\parallel(w_{tt})_{\nu_{\Lambda}}\parallel^2_{L^2(\Gamma_c)}
\end{equation*}

As in Subsection 3.2, we differentiate the full system \eqref{eq2.1}, \eqref{eq2.2} and \eqref{eq2.6}-\eqref{eq2.8} in time twice to have
\begin{eqnarray}
&&v^{i}_{ttt}-\partial_{tt}\partial_{j}(a^{jl}a^{kl}\partial_{k}v^i)+\partial_{tt}\partial_{k}(a^{ki}q)=0 \quad\mbox{in}\quad \Omega_f\times (0,T),\label{eq4.44}\\
&&a_{i}^{k}\partial_{k}v^i_{tt}+2\partial_ta^{ki}\partial_{k}v^i_t+\partial_{tt}a^{ki}\partial_{k}v^i=0  \quad\mbox{in}\quad \Omega_f\times (0,T), \quad i,j=1,2,3 \label{eq4.45}\\  \nonumber
\end{eqnarray}
and
\begin{equation}\label{eq4.46}
w_{tttt}^{i}-\div G(x)\na w^i_{tt}+w^i_{tt}=0\quad\mbox{in}\quad \Omega_e\times (0,T),\quad i=1,2,3.
\end{equation}
Moreover, the boundary conditions to \eqref{eq4.44}-\eqref{eq4.46} are changed into \begin{eqnarray}
&&w_{ttt}=v_{tt}-\gamma(w_{tt})_{\nu_\Lambda},\quad\mbox{on}\quad \Gamma_c\times(0,T),\label{eq4.47}\\
&&(w_{tt})_{\nu_\Lambda}^i=\partial_{tt}(a^{jl}a^{kl}\partial_{k}v^{i})\nu_j-\partial_{tt}(a^{ki}q)\nu_k,
\quad\mbox{on}\quad\Gamma_c\times(0,T),\label{eq4.48}\\  \nonumber
\end{eqnarray}
and
\begin{equation}\label{eq4.49}
v_{tt}=0,\quad\mbox{on}\quad \Gamma_f\times(0,T).
\end{equation}

Now we begin to derive the energy estimates for the third level energy.
\begin{lem}\label{l4.8}
Under the assumptions of Theorem \ref{th2.1}, the third level energy inequality
\begin{equation}\label{eq4.50}
E_2(t)+\int_0^{t}D_2(\tau)d\tau\leqslant E_2(0)+\int^t_0(R_2(\tau),v_{tt}(\tau))d\tau,
\end{equation}
holds for $t\in[0,T]$
where
\begin{eqnarray*}
&&\int^t_0(R_2(\tau),v_{tt}(\tau))d\tau=
2\int_0^{t}\int_{\Omega_f}\partial_t(a^{jl}a^{kl})\partial_{k}v^i_t\partial_{j}v^{i}_{tt}dxd\tau
+\int_0^{t}\int_{\Omega_f}\partial_{tt}(a^{jl}a^{kl})\partial_{k}v^i\partial_{j}v^{i}_{tt}dxd\tau \\
&&-\int_0^{t}\int_{\Omega_f}\partial_{tt}(a^{ki}q)\partial_kv^{i}_{tt}dxd\tau.
\end{eqnarray*}
\end{lem}
{\bf Proof.} First, we multiply \eqref{eq4.44} with $v^i_{tt}$ and integrate over $\Omega_f$. Sum for $i=1,2,3$, and  after integrating by parts
we obtain
\begin{eqnarray}\label{eq4.51}
&&\frac{1}{2}\frac{d}{dt}\parallel v_{tt}\parallel^2_{L^2}+\int_{\Omega_f}\partial_{tt}(a^{jl}a^{kl}\partial_{k}v^i)\partial_{j}v^{i}_{tt}dx
+\int_{\Gamma_c}\langle(w_{tt})_{\nu_\Lambda},v_{tt}\rangle d\sigma \\ \nonumber
&&-\int_{\Omega_f}\partial_{tt}(a^{ki}q)\partial_kv^{i}_{tt}dx=0.
\end{eqnarray}
Then multiply \eqref{eq4.46} by $\omega^i_{ttt}$ and integrate over $\Omega_e$ to have \begin{eqnarray}\label{eq4.52}
&&\frac{1}{2}\frac{d}{dt}(\parallel w_{ttt}\parallel^2_{L^2}+\parallel w_{tt}\parallel^2_{L^2}+\int_{\Omega_e}\<G\na w_{tt},\na w_{tt}\>dx) \\ \nonumber
&&-\int_{\Gamma_c}\langle(w_{tt})_{\nu_\Lambda},v_{tt}-\gamma(w_{tt})_{\nu_\Lambda}\rangle d\sigma=0
\end{eqnarray}
Adding \eqref{eq4.51} and \eqref{eq4.52} together leads to
\begin{eqnarray*}
&&\frac{1}{2}\frac{d}{dt}(\parallel v_{tt}\parallel^2_{L^2}+\parallel w_{ttt}\parallel^2_{L^2}+\parallel w_{tt}\parallel^2_{L^2}+\int_{\Omega_e}\<G\na w_{tt},\na w_{tt}\>dx)
+\gamma\parallel(w_{tt})_{\nu_\Lambda}\parallel^2_{L^2(\Gamma_c)}  \\
&&-\int_{\Omega_f}\partial_{tt}(a^{ki}q)\partial_kv^{i}_{tt}dx
+\int_{\Omega_f}a^{jl}a^{kl}\partial_{k}v^i_{tt}\partial_{j}v^{i}_{tt}dx
+2\int_{\Omega_f}\partial_t(a^{jl}a^{kl})\partial_{k}v^i_{t}\partial_{j}v^{i}_{tt}dx  \\
&&+\int_{\Omega_f}\partial_{tt}(a^{jl}a^{kl})\partial_{k}v^i\partial_{j}v^{i}_{tt}dx=0
\end{eqnarray*}

Finally, by the ellipticity of $\mathbf{a}(x,t)$, we integrate in time from $0$ to $t$ and obtain (\ref{eq4.50}). \hfill$\Box$\\

First we use $w_{tt}$ as multiplier as before and we get for any $t\in[0,T]$,
\begin{eqnarray}\label{eq4.53}
&&\int_0^t\int_{\Omega_e}\<G\na w_{tt},\na w_{tt}\>dxd\tau+\int_0^t\parallel w_{tt}\parallel^2_{L^2}d\tau
\leqslant
E_2(t)+E_2(0) \\ \nonumber
&&+\int_0^t\parallel w_{ttt}\parallel^2_{L^2}d\tau
+\int_0^t\langle(w_{tt})_{\nu_{\Lambda}},w_{tt}\rangle_{L^2(\Gamma_c)}d\tau.
\end{eqnarray}

Then, we continue carrying out geometric multiplier method with $H(w^i_{tt})$ and $\div H w_{tt}^i$ as multipliers.
We  obtain the following.
\begin{lem}\label{l3.9}
For every $t\in[0,T]$, the following estimate holds under the assumptions in Theorem $\ref{th2.1},$ that is
\begin{equation}\label{eq4.54}
E_2(t)+\int_0^tE_2(\tau)d\tau\leqslant
C(E_2(0)+\int^t_0(R_2(\tau),v_{tt}(\tau))d\tau).
\end{equation}
\end{lem}

\begin{rem}
We can easily generalize the estimates in Lemma \ref{l4.1}-\ref{l3.9} to any time step $[\sigma,t],$ where $0\leqslant\sigma<t\leqslant T.$
\end{rem}

\subsection{Superlinear estimates}
\hskip\parindent
Our aim of this subsection is to deal with the perturbation terms generated in the second and third level energy estimates. They are
\begin{equation*}
\int^t_0(R_1(\tau),v_t(\tau))d\tau
\end{equation*}
and
\begin{equation*}
\int^t_0(R_2(\tau),v_{tt}(\tau))d\tau.
\end{equation*}
Here we only list the results. Their proofs are similar to those in \cite{IKLT2}.

\begin{lem}$\cite[Lemma \,4.10]{IKLT2}$\label{l3.10}
We have
\begin{equation}
\mid(R_1(t),v_t)\mid\leqslant C\parallel v\parallel_{H^1}^{\frac{1}{2}}\parallel v\parallel_{H^2}^{\frac{1}{2}}\parallel v_t\parallel_{H^1}(\parallel v\parallel_{H^2}+\parallel q\parallel_{H^1})+C\parallel v\parallel_{H^1}^{\frac{3}{2}}\parallel v\parallel_{H^2}^{\frac{1}{2}}\parallel q_t\parallel_{H^1},
\end{equation}
for all $t\in[0,T]$.
\end{lem}

\begin{lem}$\cite[Lemma\, 4.11]{IKLT2}$\label{l3.11}
For $\epsilon_0\in(0,\frac{1}{C}]$, we have
\begin{eqnarray*}
\int_0^t(R_2(s),v_{tt}(s))ds&\leqslant& \epsilon_0\int_0^t\parallel \nabla v_{tt}\parallel_{L^2}^{2}ds
+C_{\epsilon_0}\int_0^t(\parallel v\parallel_{H^3}^{2}+\parallel q\parallel_{H^2}^{2})(\parallel v\parallel_{H^1}^{\frac{5}{2}}\parallel v\parallel_{H^3}^{\frac{3}{2}}+\parallel v_t\parallel_{H^1}^2)ds \\ \nonumber
&+&C_{\epsilon_0}\int_0^t\parallel v\parallel_{H^1}^{\frac{3}{2}}\parallel v\parallel_{H^3}^{\frac{1}{2}}\parallel q_t\parallel_{H^1}^2ds+\epsilon_0\parallel q_t(t)\parallel_{H^1}^2+\epsilon_0\parallel v_t(t)\parallel_{H^2}^2 \\ \nonumber
&+&\epsilon_0\parallel v(t)\parallel_{H^3}^2+C_{\epsilon_0}\parallel v(t)\parallel_{H^1}^{6}\parallel v(t)\parallel_{H^2}^{4}+C_{\epsilon_0}\parallel v(t)\parallel_{H^1}^{2}\parallel v(t)\parallel_{H^2}^{2}\parallel v_t(t)\parallel_{L^2}^2   \\ \nonumber
&+&C\int_0^t(\parallel v\parallel_{H^2}^2+\parallel v_t\parallel_{H^1}^{\frac{1}{2}}\parallel v_t\parallel_{H^2}^{\frac{1}{2}})\parallel q_t\parallel_{H^1}\parallel v_t\parallel_{H^1}ds   \\ \nonumber
&+&C\int_0^t(\parallel v\parallel_{H^2}^3+\parallel v_t\parallel_{H^1}\parallel v\parallel_{H^1}^{\frac{1}{4}}\parallel v\parallel_{H^3}^{\frac{3}{4}})\parallel q_t\parallel_{H^1}\parallel v\parallel_{H^1}^{\frac{3}{4}}\parallel v\parallel_{H^3}^{\frac{1}{4}}ds   \\ \nonumber
&+&C\parallel v(0)\parallel_{H^3}^{6}+C\parallel v_t(0)\parallel_{H^1}^{4}+C\parallel q_t(0)\parallel_{H^1}^2,
\end{eqnarray*}
for all $t\in[0,T]$.
\end{lem}

\begin{rem}
From the proof in \cite{IKLT2} the superlinear estimate in Lemma \ref{l3.11} holds true to any time step $[\sigma,t],$ where $0\leqslant\sigma<t\leqslant T.$
\end{rem}

\section{Energy decay and global existence of the system}
\setcounter{equation}{0}
\hskip\parindent
We focus on the global existence of solutions and the energy decay estimates.

Let
\begin{equation*}
\mathcal{X}(t)=E(t)+E_1(t)+E_2(t)+\hat{\epsilon}_1(\parallel\nabla v\parallel^2_{L^2}+\parallel\nabla v_t\parallel^2_{L^2}),
\end{equation*}
where $\hat{\epsilon}_1>0$  be given sufficiently small to be determined later.

We make some preparations for the proof of Theorem \ref{th2.1} as that in \cite{IKLT2}.

We have
\beq
\parallel\nabla v(t)\parallel^2_{L^2}&&=\parallel\nabla v(0)\parallel^2_{L^2}+\int^t_0\frac{d}{d\tau}\parallel\nabla v(\tau)\parallel^2_{L^2}d\tau\nonumber\\
&&= \parallel\nabla v(0)\parallel^2_{L^2}+2\int^t_0\parallel\nabla v\parallel_{L^2}\parallel\nabla v_t\parallel_{L^2}d\tau\nonumber\\
&&\leq\parallel\nabla v(0)\parallel^2_{L^2}+C\int^t_0(D(\tau)+D_1(\tau))d\tau,    \label{eq5.3}
\eeq
Similarly, we obtain
\begin{equation}\label{eq5.4}
\parallel\nabla v_t(t)\parallel^2_{L^2}\leqslant\parallel\nabla v_t(0)\parallel^2_{L^2}+C\int^t_0(D_1(\tau)+D_2(\tau))d\tau.
\end{equation}

In addition, it follows from Lemmas \ref{l4.1} and \ref{l4.3} that
\begin{equation}\label{eq5.5}
E(t)+\int^t_0E(\tau)d\tau+\int^t_0D(\tau)d\tau\leqslant CE(0).
\end{equation}
Combining Lemmas \ref{l4.5}, \ref{l4.7} and \ref{l3.10}, we have
\begin{equation}\label{eq5.6}
E_1(t)+\int^t_0E_1(\tau)d\tau+\int^t_0D_1(\tau)d\tau\leqslant CE_1(0)+\int^t_0\mathbf{P}_1(\parallel v\parallel_{H^2},\parallel q\parallel_{H^1},\parallel v_t\parallel_{H^1},\parallel q_t\parallel_{H^1})d\tau,
\end{equation}
and, from Lemmas \ref{l4.8}, \ref{l3.9} and \ref{l3.11},
\begin{eqnarray}\label{eq5.7}
E_2(t)&+&\int^t_0E_2(\tau)d\tau+\int^t_0D_2(\tau)d\tau\leqslant CE_2(0)+\epsilon_0\parallel v(t)\parallel_{H^3}^2
+\epsilon_0\parallel q_t(t)\parallel_{H^1}^2  \nonumber\\
&+&\epsilon_0\parallel v_t(t)\parallel_{H^2}^2
+\epsilon_0\int_0^t\parallel \nabla v_{tt}\parallel_{L^2}^{2}d\tau
+\mathbf{P}_2(\parallel v\parallel_{H^2},\parallel v_t\parallel_{L^2})  \nonumber\\
&+&\int^t_0\mathbf{P}_3(\parallel v\parallel_{H^3},\parallel q\parallel_{H^2},\parallel v_t\parallel_{H^2},\parallel q_t\parallel_{H^1})d\tau   \nonumber\\
&+&\mathbf{P}_4(\parallel v(0)\parallel_{H^3},\parallel v_t(0)\parallel_{H^1},\parallel q_t(0)\parallel_{H^1}),
\end{eqnarray}
where the symbols $\mathbf{P}_1,\mathbf{P}_2,\mathbf{P}_3$ and $\mathbf{P}_4$  denote the superlinear polynomials of
their arguments, which are allowed to depend on $\epsilon_0$ from Lemmas \ref{l3.10} and \ref{l3.11}. Now multiply \eqref{eq5.3} and \eqref{eq5.4} by
$\hat{\epsilon}_1$, sum  \eqref{eq5.5}-\eqref{eq5.7} and then add them together to obtain
\begin{eqnarray}\label{eq5.8}
&&\mathcal{X}(t)+\int^t_0\mathcal{X}(\tau)d\tau\leqslant C\mathcal{X}(0)+\epsilon_0\parallel v(t)\parallel_{H^3}^2
+\epsilon_0\parallel q_t(t)\parallel_{H^1}^2+\epsilon_0\parallel v_t(t)\parallel_{H^2}^2 \nonumber\\
&&+\mathbf{P}_1(\parallel v\parallel_{H^2},\parallel v_t\parallel_{L^2})
+\int^t_0\mathbf{P}_2(\parallel v\parallel_{H^3},\parallel q\parallel_{H^2},\parallel v_t\parallel_{H^2},\parallel q_t\parallel_{H^1})d\tau  \nonumber\\
&&+\mathbf{P}_3(\parallel v(0)\parallel_{H^3},\parallel v_t(0)\parallel_{H^1},\parallel q_t(0)\parallel_{H^1}).
\end{eqnarray}

From \eqref{eq3.9}, we have
\begin{equation}\label{eq5.9}
\parallel v\parallel_{H^2}^2+\parallel q\parallel_{H^1}^2\leqslant C\mathcal{X}(t).
\end{equation}
Thanks to \eqref{eq3.8} and \eqref{eq5.9}, we obtain
\begin{equation}\label{eq5.10}
\parallel v\parallel_{H^3}^2+\parallel q\parallel_{H^2}^2\leqslant C\mathcal{X}(t).
\end{equation}

From \eqref{eq3.10} and \eqref{eq5.10}, we deduce that
\begin{equation}\label{eq5.11}
\parallel v_t\parallel_{H^2}^2+\parallel q_t\parallel_{H^1}^2\leqslant C\mathcal{X}(t)+C\mathcal{X}^2(t).
\end{equation}

Using  \eqref{eq5.10} and \eqref{eq5.11} and choosing $\epsilon_0$ sufficiently small, we  deduce from \eqref{eq5.8}
that
\begin{equation}\label{eq5.12}
\mathcal{X}(t)+\int^t_0\mathcal{X}(\tau)d\tau\leqslant C\mathcal{X}(0)+\mathbf{P}(\mathcal{X}(t))+\int^t_0\mathbf{P}(\mathcal{X}(\tau))d\tau+\mathbf{P}(\mathcal{X}(0)),
\end{equation}
where $\mathbf{P}$ is a superlinear polynomial as well. We  rewrite \eqref{eq5.12} as
\begin{equation}\label{eq5.13}
\mathcal{X}(t)+\int^t_{0}\mathcal{X}(\tau)d\tau\leqslant C_0\sum_{j=1}^{m}\int^t_{0}\mathcal{X}(\tau)^{\alpha_j}d\tau+C_0\sum_{k=1}^{n}\mathcal{X}(t)^{\beta_k}
+C_0\sum_{k=1}^{n}\mathcal{X}(0)^{\beta_k}+C_0\mathcal{X}(0),
\end{equation}
where $C_0\geqslant 1$, $\alpha_1,...,\alpha_m>1$ and $\beta_1,...,\beta_n>1$.

Following \cite[Lemma\, 5.1]{IKLT2}, we have
\begin{lem}\label{l5.1}
Suppose that $\mathcal{X}:[0,\infty)\rightarrow[0,\infty)$ is continuous for all $t$ such that $\mathcal{X}(t)$ is finite and assume that it satisfies
\begin{equation*}
\mathcal{X}(t)+\int^t_{\tau}\mathcal{X}(s)ds\leqslant C_0\sum_{j=1}^{m}\int^t_{\tau}\mathcal{X}(s)^{\alpha_j}ds+C_0\sum_{k=1}^{n}\mathcal{X}(t)^{\beta_k}
+C_0\sum_{k=1}^{n}\mathcal{X}(\tau)^{\beta_k}+C_0\mathcal{X}(\tau),
\end{equation*}
where $\alpha_1,...,\alpha_m>1$ and $\beta_1,...,\beta_n>1.$ Also, assume that $\mathcal{X}(0)\leqslant\epsilon$. If $\epsilon\leqslant\frac{1}{C}$, where the constant $C$ depends on $C_0, m, \alpha_1,...,\alpha_m, \beta_1,...,\beta_n$, we have
$\mathcal{X}(t)\leqslant C\epsilon e^{-\frac{t}{C}}$.
\end{lem}

{\bf Proof of Theorem \ref{th2.1}}\,\,\, Using Lemma \ref{l5.1} and following the proof of \cite[Theorem 2.1]{IKLT2}, we complete the proof of Theorem \ref{th2.1}.\hfill$\Box$\\


\begin{thebibliography}{99}
\bibitem{CS1} Coutand. D, Shkoller. S, Motion of an elastic solid inside an incompressible viscous fluid, Arch. Ration. Mech. Anal. \textbf{176}(2005), 25-102.
\bibitem{CS2} Coutand. D, Shkoller. S, The interaction between quasilinear elastodynamics and the Navier-Stokes equations, Arch. Ration. Mech. Anal. \textbf{179}(2006), 303-52.
\bibitem{feng2001}Feng S., Feng D., Boundary stabilization of wave equations with variable coefficients, Science in China.
    Vol.44 No.3 (2001)
\bibitem{IKLT1} Ignatova M, Kukavica I, Lasiecka I and Tuffaha A, On well-posedness for a free boundary fluid-structure model, J. Math. Phys. \textbf{53}(2012), 115624
\bibitem{IKLT2} Ignatova M, Kukavica I, Lasiecka I and Tuffaha A, On well-posedness and small data global existence for an interface damped free boundary fluid-structure model, {\it Nonlinearity,} \textbf{27}(2014), 467-99.
\bibitem{IKLT3} Ignatova M, Kukavica I, Lasiecka I, and Tuffaha A, Small data global existence for a fluid-structure model,{\it  Nonlinearity}, \textbf{30}(2017), 848-98.
\bibitem{KT1} Kukavica I, Tuffaha A, Solutions to a fluid-structure interaction free boundary problem, Discrete. Contin. Dyn. Syst. \textbf{32}(2012), 1355-89.
\bibitem{KT2} Kukavica I, Tuffaha A, Regularity of solutions to a free boundary problem of fluid-structure interaction, {\it Indina. Univ. Math. J.} \textbf{61}(2012), 1817-59.
\bibitem{LTY}  Lasiecka I, Triggiani R., Yao P. F., Inverse/observability estimates for second-order hyperbolic equations with variable coefficients. {\it J. Math. Anal. Appl.} 235 (1999), no. 1, 13-57.
\bibitem{TY}  Triggiani R, Yao P. F., Carleman estimates with no lower-order terms for general Riemann wave equations. Global uniqueness and observability in one shot. Special issue dedicated to the memory of Jacques-Louis Lions. {\it Appl. Math. Optim.} 46 (2002), no. 2-3, 331-375.
\bibitem{Yao 99} Yao P. F.,  On the observability inequalities for exact controllability of wave equations with variable coefficients. SIAM {\it J. Control Optim.} 37 (1999), no. 5, 1568-1599.
\bibitem{Yao 2011} Yao P. F., {\it  Modeling and Control in Vibrational and Structural Dynamics. A Differential Geometric Approach},  CRC Press, Boca Raton, Florida, 2011.
\bibitem{Yao 2007}Yao P. F., Global smooth solutions for the quasilinear wave equation with boundary dissipation, {\it J. Diff. Eq.} 241 (2007),  62-93.
\bibitem{wh 1989} Wu H., Shen C. L.,  Yu Y. L., {\it An Introduction to Riemannian Geometry (in Chinese)}, Beijing University Press, Beijing, 1989.
\bibitem{dCa} Carmo M. Do., {\it Riemannian Geometry}, Birkha\"{u}ser, Basel, 1992.
\bibitem{GHL} Gallo S., Hulin D.,  Lafontaine J.,  {\it Riemannian Geometry}, Springer, 2004.
\end{thebibliography}
\end{document}